\documentclass[12pt]{article}
\usepackage[a4paper, total={7in, 10in}]{geometry}
\usepackage[parfill]{parskip}
\usepackage{physics, tensor, float, subcaption}
\usepackage{graphicx}
\graphicspath{ {Plots/} }
\usepackage{jhep-mod}
\usepackage{bm}
\usepackage{soul}
\usepackage{amssymb,amsmath,amsthm}
\usepackage{mathrsfs}
\usepackage[utf8]{inputenc}
\usepackage{enumerate}
\usepackage{bigints}
\usepackage{xcolor}
\usepackage{appendix}
\usepackage{graphicx}
\usepackage{float}
\usepackage{tikz}
\usepackage{setspace}
\usepackage{cancel}
\usepackage{doi}
\definecolor{purple}{rgb}{1,0,1}
\definecolor{lime}{HTML}{A6CE39} 

\newcommand{\blue}[1]{{\color{blue} #1}}

\definecolor{lime}{HTML}{A6CE39}
\newcommand{\orcidicon}{%
	\begin{tikzpicture}
	\draw[lime, fill=lime] (0,0) 
		circle [radius=0.16] 
		node[white] {{\fontfamily{qag}\selectfont \tiny ID}};
	\draw[white, fill=white] (-0.0625,0.095) 
		circle [radius=0.007];
	\end{tikzpicture}
	\hspace{-5mm}
}
\newcommand\orcidMatt{{\href{https://orcid.org/0000-0003-1088-6485}{\orcidicon}}}

\renewcommand{\O}{\mathcal{O}}
\begin{document}

\title{
\leftline{\blue{On the arithmetic average of the first $n$ primes}}
}

\author{
\Large
Matt Visser\!\orcidMatt\!
}
\affiliation{School of Mathematics and Statistics, Victoria University of Wellington, \\
\null\qquad PO Box 600, Wellington 6140, New Zealand.}
\emailAdd{matt.visser@sms.vuw.ac.nz}
\def\theta{\vartheta}
\def\O{{\mathcal{O}}}
\def\Li{{\mathrm{Li}}}

\abstract{
\vspace{1em}

The arithmetic average of the first $n$ primes, $\bar p_n = {1\over n} \sum_{i=1}^n p_i$, exhibits very many \break interesting and subtle properties. Since the transformation from $p_n \to \bar p_n$ is extremely easy to invert, $p_n = n\bar p_n - (n-1)\bar p_{n-1}$, it is clear that these two sequences $p_n \longleftrightarrow \bar p_n$ must ultimately carry exactly the same information. But the averaged sequence $\bar p_n$, while very closely correlated with the primes, ($\bar p_n \sim {1\over2} p_n$),  is much ``smoother'', and much better behaved. 
Using extensions of various standard results I shall demonstrate that the  prime-averaged sequence $\bar p_n$ satisfies prime-averaged analogues of the Cramer, Andrica, Legendre, Oppermann, Brocard, Fourges, Firoozbakht, Nicholson, and Farhadian conjectures. (So these prime-averaged analogues are not conjectures, they are theorems.) The crucial key to enabling this pleasant behaviour is the ``smoothing'' process inherent in averaging. Whereas the asymptotic behaviour of the two sequences is very closely correlated the local fluctuations are quite different. 

\bigskip

\bigskip
\noindent
{\sc Date:} 16 July 2025; \LaTeX-ed \today


\bigskip
\noindent
{\sc arXiv:} {\href{https://arxiv.org/abs/2505.04951}{\color{blue}2505.04951}}     [math.NT]

\bigskip
\noindent
{\sc Published:}  Mathematics {\bf 13} (2025) 2279.
\doi{10.3390/math13142279}. 
}

\maketitle
\def\tr{{\mathrm{tr}}}
\def\diag{{\mathrm{diag}}}
\def\cof{{\mathrm{cof}}}
\def\pdet{{\mathrm{pdet}}}
\def\d{{\mathrm{d}}}
\parindent0pt
\parskip7pt
\def\Kerr{{\scriptscriptstyle{\mathrm{Kerr}}}}
\def\eos{{\scriptscriptstyle{\mathrm{eos}}}}

\clearpage
\section{Introduction and background}
\def\half{{\textstyle{1\over2}}}
\def\third{{\textstyle{1\over3}}}
\def\quarter{{\textstyle{1\over4}}}
\def\fifth{{\textstyle{1\over5}}}
\def\ninth{{\textstyle{1\over9}}}
\def\Heaviside{{\mathrm{Heaviside}}}
\counterwithout{equation}{section}

Consider the average of the first $n$ primes
\begin{equation}
\bar p_n =  {1\over n} \sum_{i=1}^n p_i.
\end{equation}
Explicitly
\begin{equation}
\bar p_n \in \left\{2,  2\half, 3\third, 4\quarter, 5{{\textstyle{3\over5}}}, 6{{\textstyle{5\over6}}}, 
8{{\textstyle{2\over7}}}, 9{{\textstyle{5\over8}}}, 11{{\textstyle{1\over9}}}, 
12{{\textstyle{9\over10}}},\dots
\right\} \hbox{ for } n\in\{1,2,3,\dots\}.
\end{equation}
Quite a lot is already known about this sequence~\cite{A034387,Dusart:1998-thesis,Axler, Axler2, Hassani,Sinha,Mandl,A351914,Sun}.

For instance the asymptotic behaviour is known to be (\hbox{OEIS A034387~\cite{A034387}})
\begin{equation}
\bar p_n \sim {p_n\over 2} \sim {n \ln n\over 2}.
\end{equation}
More precisely (\hbox{Dusart~\cite{Dusart:1998-thesis}})
\begin{equation}
\bar p_n = {1\over2} \left( n\ln(n\ln n) - {3n\over2} + o(1) \right).
\end{equation}
Several higher-order terms in the asymptotic expansion are also known~\cite{Axler, Axler2, Hassani, Dusart:1998-thesis, Sinha}, 
but are not of direct relevance to the present discussion.

Various other explicit upper and lower bounds are also known
\begin{eqnarray}
&\bar p_n  < {p_n\over2}; \qquad\qquad 
&(p_n> 19; n \geq 9; \;\hbox{ Mandl~\cite{Mandl}, OEIS A351914~\cite{A351914}});
\\
&\bar p_n  < {p_n\over2} - {n\over 14}; \qquad\qquad 
&(n\geq 10; \;\hbox{ Hassani~\cite{Hassani}});
\\
&\bar p_n > p_{[n/2]}; \qquad\qquad 
&(n\geq 2; \; \hbox{ Dusart~\cite{Dusart:1998-thesis}}).
\end{eqnarray}
Herein, I wish to take the discussion in a rather different direction. 

Since $\bar p_n \sim {p_n\over 2}$ these two sequences are very closely correlated.
Furthermore, since the relation, $p_n \longleftrightarrow \bar p_n$ is easily invertible $p_n = n\bar p_n - (n-1)\bar p_{n-1}$, these two sequences ultimately carry identical information. 
Nevertheless, as discussed below, the two sequences also exhibit profound differences.

\enlargethispage{50pt}
Specifically, I shall demonstrate below that the  prime-averaged sequence $\bar p_n$ \emph{satisfies} suitably formulated prime-averaged analogues of the Cramer, Andrica, Legendre, Oppermann, Brocard, Fourges, Firoozbakht, Nicholson, and Farhadian
 conjectures. (So these prime-averaged analogues are not conjectures, they are theorems.) 
The crucial observation enabling this pleasant behaviour is the ``smoothing'' process inherent in averaging, we shall soon se that the prime-averaged gaps are in a suitable sense extremely small. On the other hand, while the asymptotic behaviour of the two sequences is very tightly correlated, we shall see that  the fluctuations are quite different. 

\clearpage
\section{Standard and easy results}

In this section I will introduce  a few basic tools (\hbox{Rosser~\cite{Rosser1}, \hbox{Rosser--Schoenfeld~\cite{Rosser+Schoenfeld}}}):
\begin{eqnarray}
&p_n > n\ln n;  \qquad &(n\geq1);
\\
&p_n < n\ln(n\ln n);  \qquad\qquad &(n\geq6).
\end{eqnarray}
It is sometimes useful to note 
\begin{equation}
\ln(n\ln n) = \ln n + \ln\ln n = \ln n \left(1 + {\ln\ln n\over\ln n} \right) \leq (1+e^{-1}) \; \ln n.
\end{equation}
Consequently 
\begin{equation}
p_n < (1+e^{-1}) \; n\ln n;  \qquad (n\geq 4).
\end{equation}

{\bf \underline{Lemma:}} Still easy, (but perhaps somewhat less well known), is the first slightly non-trivial (post Bertrand--Chebyshev) bound on the $n^{th}$ prime gap:
\begin{equation}
g_n = p_{n+1}-p_n < n; \qquad \hbox{so that} \qquad p_{n+1} < p_n+ n; \qquad (n\geq 1).
\end{equation}
{\bf \underline{Proof:}} Rosser and Schoenfeld~\cite{Rosser+Schoenfeld} give $|\theta(x)-x| <{1\over2}\; {x\over\ln x} $ for $x\geq 563=p_{103}$.  Thence, (checking small values of $x$ by explicit computation), $|\theta(x)-x| <{1\over2}\; \pi(x) $ for $x\geq 347=p_{69}$. But then $|\theta(p_n)-p_n|<{1\over2} n$ for $n\geq 26$.
By evaluating $|\theta(x)-x|$ at $x=p_{n+1}^-$, just below the $(n+1)^{th}$ prime, we also see $|\theta(p_n)-p_{n+1}|<{1\over2} n$ for $n\geq 69$. Thence by the triangle inequality $p_{n+1}-p_n< n$, certainly for all integers $n\geq 69$. Checking lower integers by explicit computation, $p_{n+1}-p_n< n$ for $n\geq 1$. 
\hfill $\Box$

Note that this in turn implies 
\begin{equation}
p_{n+1} < p_n +n < n [ \ln(n\ln n) +1 ]; \qquad (n\geq 5). 
\end{equation}

For another lower bound on $\bar p_n$ note
\begin{equation}
\bar p_n > {1\over n} \sum_{i=1}^n i \ln i =
{1\over n} \sum_{i=2}^n i \ln i.
\end{equation}
We now bound this sum by the integral
\begin{equation}
\bar p_n  > {1\over n} \int_1^n u \ln u \;d u = {n\ln n\over2} - {n\over4}+{1\over4n}
= {n\over2} \left( \ln n -{1\over2}\right) +{1\over4n},
\end{equation}
with this inequality holding for $n\geq 1$.

Finally, we shall also have occasion to use (\hbox{Dusart~\cite{Dusart:2018}})
\begin{equation}
\pi(x) > {x\over\ln x -1}; \qquad (x>5393=p_{711}),
\end{equation}
and (\hbox{Dusart~\cite{Dusart:2018}})
\begin{equation}
\pi(x) < {x\over\ln x -{139\over125}}; \qquad (x>e^{139/125}).
\end{equation}

\section{Counting the averaged primes}
Define
\begin{equation}
\bar\pi(x):= \#\{ i: \bar p_i \leq x\}.
\end{equation}
This counts the averaged-primes $\bar p_n$. (I emphasize that here we are counting the averaged-primes, not averaging the count of primes. That would instead be something akin to $\pi_\mathrm{average}(x) = {1\over \pi(x)} \sum_{i=1}^{\pi(x)} \pi(i)$, a rather different quantity.)

\begin{itemize}
\item 
From $\bar p_n  < {p_n\over2}$ we see $\{ i: \bar p_i \leq x\}  \supseteq \{ i: {p_i\over2}  \leq x\} =  \{ i: p_i \leq 2x\}$, and so for $x \geq 19$ we have
 $\#\{ i: \bar p_i \leq x\}  \geq \#\{ i: p_i \leq 2x\}$. 
 In terms of the average-prime counting function this implies 
\begin{equation}
\bar\pi(x) \geq \pi(2x); \qquad \qquad (x> 19=p_8).
\end{equation}
Explicitly checking smaller values of $x$ this saturates the domain of validity.

\item From $\bar p_n > p_{[n/2]}$ we see  $\{ i: \bar p_i \leq x\}  \subseteq \{ i: p_{[i/2]} \leq x\} $ whence for $n\geq 2$, (implying $x \geq 3$), we have:
 \begin{equation}
\bar\pi(x) \leq \#\{ i: p_{[i/2]} \leq x\} \leq  2 \; \#\{ i: p_{i} \leq x\}\leq 2\; \pi(x). 
\end{equation}
Explicitly checking smaller values of $x$ this saturates the domain of validity.
\begin{equation}
\bar\pi(x) \leq 2\pi(x); \qquad \qquad (x\geq 3=p_2).
\end{equation}
\end{itemize}

Combining, we see that we have reasonably tight upper and lower bounds on the count of averaged primes:
 \begin{equation}
(x> 19=p_8) \qquad\qquad\pi(2x) \leq \bar\pi(x) \leq 2 \pi(x) \qquad \qquad (x\geq 3=p_2).
\end{equation}

\section{Bounding the gaps in the averaged primes}

\enlargethispage{40pt}
First note that
\begin{equation}
\label{E:bertrand}
\bar p_{n+1} = {\sum_{i=1}^{n+1} p_i \over n+1} = {\sum_{i=1}^{n} p_i + p_{n+1} \over n+1} = {n \bar p_n + p_{n+1}\over n+1} = \bar p_n + {p_{n+1}-\bar p_n\over n+1}.
\end{equation}
In view of Mandl's inequality ($\bar p_n < p_n/2$ for $n\geq 9$, see~\cite{Mandl,A351914}) that last term is definitely positive and the sequence $\bar p_n$ is monotone increasing. 
Let us define the gap in the averages as
\begin{equation}
\bar g_n =\bar p_{n+1} - \bar p_n.
\end{equation} 
While the average primes $\bar p_n$ are monotone increasing, the gaps in the averages $\bar g_n$ have no nice monotonicity properties.

(I emphasise, $\bar g_n$ is the gap in the averages, not the averages of the gaps --- that would instead be 
\begin{equation}
(g_n)_\mathrm{average} = {1\over n} \sum_{i=1}^n g_i = {1\over n} \sum_{i=1}^n (p_{i+1}-p_i) 
= {1\over n} (p_{n+1}-p_1)  = {1\over n} (p_{n+1}-2),
\end{equation}
a quantity interesting in its own right, but for totally different reasons.)

For the gap in the averages we have
\begin{equation}
\bar g_n = \bar p_{n+1} - \bar p_n = {p_{n+1}-\bar p_n\over n+1} < {n [ \ln(n\ln n) +1 ] - {n\over2} \left( \ln n -{1\over2}\right)\over n+1},
\end{equation}
with the inequality holding for $n\geq 5$.  Rearranging
\begin{equation}
\bar g_n  <  \ln n \left\{ {n \over n+1}  \;\; 
\left[ \left(1 +{\ln \ln n\over n} + {1\over\ln n}\right) - {1\over2} \left(1-{1\over 2\ln n}\right) \right] \right\},
\end{equation}
that is
\begin{equation}
\bar g_n  <  {\ln n\over 2} \left\{ {n \over n+1}  \;\; 
\left[1 +{2\ln \ln n\over n} + {5\over2\ln n} \right] \right\}, \qquad (n\geq 5).
\end{equation}
Within the domain of validity of this inequality, the quantity in braces is monotone decreasing (but always exceeds unity) and, searching for a value of $n$ such that quantity in braces drops below  2,  we certainly have $\bar g_n  <  \ln n$ for  $n\geq 439$.
Explicitly checking the smaller integers one finds
\begin{equation}
\bar g_n  <  \ln n , \qquad (n\geq 3).
\end{equation}
In the other direction
\begin{equation}
\bar g_n = {p_{n+1} -\bar p_n\over n+1} > {p_{n} +2 -{1\over 2} p_n\over n+1}= {1\over2}{p_{n} +4\over n+1}
> {1\over2} {n\ln n+4 \over n+1} = {\ln n\over 2} \left\{1+{4\over\ln n}\over 1+{1\over n}\right\}.
\end{equation}
Thence, since the quantity in braces always exceeds unity, we see
\begin{equation}
\bar g_n  >  {\ln n\over2} , \qquad (n\geq 1).
\end{equation}
Overall we have rather good upper and lower bounds
\begin{equation}
(n\geq 1), \qquad {\ln n\over2} < \bar g_n  <  \ln n , \qquad (n\geq 3).
\end{equation}
Indeed, from this  we see $\bar g_n \sim \ln n$ and consequently ${\bar g_n\over\bar p_n} \sim {2\over n}$. This extremely rapid falloff in the size of the relative gaps in the averaged primes is ultimately the key observation underlying the computations below.

To be more explicit, for $n \geq 3$ we certainly have
\begin{equation}
{\bar g_n\over\bar p_n} < {\ln n\over {n\over2}(\ln n -1/2) +1/4} < {2.871142034\over n}.
\end{equation}
By somewhat rearranging and tightening the discussion, from (\ref{E:bertrand}) we have
\begin{equation}
n \; {\bar g_n\over \bar p_n} = {n \over n+1} \left({p_{n+1}\over\bar p_n}-1\right)
< {n \over n+1} \left({[n+1]\ln([n+1]\ln[n+1])\over {n\over2} \left( \ln n -{1\over2}\right) +{1\over4}}-1\right)
\end{equation}
Now this last quantity is less than 2 for $n\geq100$. Thence certainly
\begin{equation}
{\bar g_n\over \bar p_n}  < {2\over n}; \qquad (n>100).
\end{equation}
Explicitly checking smaller integers we have the stronger result
\begin{equation}
{\bar g_n\over \bar p_n}  < {2\over n}; \qquad (n\geq 4).
\end{equation}

In the other direction,  for $n\geq 6$ we have
\begin{equation}
{\bar g_n\over\bar p_n} > {{1\over2} \ln n\over {1\over 2} p_n } >{\ln n\over n \ln(n\ln n)}
= {1\over n (1+{\ln\ln n\over \ln n} )} 
\end{equation}
which is bounded below by 
\begin{equation}
{1\over n (1+e^{-1})} > {0.7310585785\over n};   \qquad (n\geq 6).
\end{equation}
Explicitly checking smaller integers
\begin{equation}
{\bar g_n\over\bar p_n} > {1\over n (1+e^{-1})} > {0.7310585785\over n};   \qquad (n\geq 3).
\end{equation}
A slightly tidier summary of these results is
\begin{equation}
(n\geq 4) \qquad\qquad {7\over 10 n} < {\bar g_n\over \bar p_n}  < {2\over n}; \qquad\qquad (n\geq 4).
\end{equation}
We shall re-use these bounds several times in the discussion below.

\section{Prime-average analogues of some standard conjectures}

Given the rather tight bound on the prime-average gaps derived above, it is now in very many cases relatively easy to formulate \emph{provable} prime-average analogues of various standard conjectures. (Often the only major difficulty lies in designing and formulating a suitable analogue.) 
Below I present  a few examples.

\subsection{Prime-average analogue of Cramer}

The ordinary Cramer conjecture is the hypothesis that the ordinary prime gaps $g_n$ satisfy $g_n = \O([\ln p_n]^2)$. 
The prime-average analogue would be $\bar g_n = \O([\ln \bar p_n]^2) $.
But since $\bar p_n = \O(n\ln n)$ we have $\ln \bar p_n = \O(\ln n)$ and 
$\O([\ln \bar p_n]^2) = \O([\ln n]^2)$.  

So the prime-average analogue of Cramer would be tantamount to making the claim $\bar g_n = \O([\ln n]^2) $.
But since we have already proved the very much stronger result that $\bar g_n = \O(\ln n) $, this is a triviality.

\subsection{Prime-average analogue of Andrica}

The ordinary Andrica conjecture is the hypothesis that the ordinary primes $p_n$ satisfy $\sqrt{p_{n+1}}- \sqrt{p_n} < 1$. The prime-average analogue would be that the averaged primes $\bar p_n$ satisfy $\sqrt{\bar p_{n+1}}- \sqrt{\bar p_n} < 1$. But this is easily checked to be true and in fact much more can be said.

Using our previous results we compute (for $n> 1$) 
\begin{equation}
\sqrt{\bar p_{n+1}} - \sqrt{\bar p_n} = { \bar p_{n+1}- \bar p_n\over \sqrt{\bar p_{n+1}} + \sqrt{\bar p_n} } 
< { \ln n \over 2 \sqrt{\bar p_n} } < { \ln n \over 2 \sqrt{{n\over2}(\ln n -{1\over2} ) } }.
\end{equation}
That is
\begin{equation}
\sqrt{\bar p_{n+1}} - \sqrt{\bar p_n} 
< \sqrt{\ln n\over n} \left\{2\left(1-{1\over2\ln n}\right)\right\}^{-1/2}.
\end{equation}
But for $n> e$ (implying $n \geq 3$) the quantity in braces lies in the range $(1,2)$ and so in this range 
\begin{equation}
\sqrt{\bar p_{n+1}} - \sqrt{\bar p_n}  < \sqrt{\ln n\over n}; \qquad (n\geq 3).
\end{equation}
Explicitly checking $n\in\{1,2\}$ in fact we see
\begin{equation}
\sqrt{\bar p_{n+1}} - \sqrt{\bar p_n}  < \sqrt{\ln n\over n}; \qquad (n\geq 2).
\end{equation}

Note this is asymptotically much stronger than just a constant bound, and we could replace it with  a considerably weaker statement with a slightly greater range of validity
\begin{equation}
\sqrt{\bar p_{n+1}} - \sqrt{\bar p_n}  \leq  \sqrt{\bar p_{3}} - \sqrt{\bar p_2} < {1\over4}; \qquad (n\geq 1).
\end{equation}
So the prime-average analogue of Andrica is unassailably true.

\subsection{Prime-average analogue of Legendre}

The ordinary Legendre conjecture is the hypothesis that the ordinary primes $p_n$ satisfy $\pi([m+1]^2)> \pi(m^2)$ for integer $m \geq 1$. The most naive prime-average analogue would be that the averaged primes $\bar p_n$ satisfy $\bar\pi([m+1]^2)> \bar\pi(m^2)$ for integer $m \geq 1$. But this is easily checked to be true and in fact much more can be said.

Recall
\begin{equation}
\pi(2x) \leq \bar\pi(x) \leq 2 \pi(x). 
\end{equation}
Let $\bar\pi(m^2)=n$, then by definition $\bar p_n$ is the largest prime-average below $m^2$, and $\bar p_{n+1}$ is the smallest prime-average above $m^2$. But then, using our previous results, for integer $n\geq 3$ corresponding to $m\geq \sqrt{5}$, we have
\begin{equation}
\bar p_{n+1} = \bar p_n + \bar g_n < m^2+ \ln n= m^2 +\ln \bar \pi(m^2) < m^2 +\ln[2  \pi(m^2)] < m^2 +\ln 2 + \ln\pi(m^2).
\end{equation}
Thence, since $\pi(m^2)$ is certainly (very much) less than $m^2$, and $\ln 2 < 1$,  we certainly have the (extremely crude) bound
\begin{equation}
\bar p_{n+1}  < m^2+\ln 2 + \ln(m^2) < m^2 + 1 + 2\ln m = (m+1)^2 - 2[m-\ln m] <(m+1)^2 .
\end{equation}
This is sufficient  to prove the prime-average version of Legendre for $m\geq \sqrt{5}$, though it is clear that the argument can be very considerably tightened. 
In particular, it is easy to check that
\begin{equation}
\bar p_n^2 \in\left\{4, 6\quarter, 11\ninth, 19{{\textstyle{1\over16}}},
31{{\textstyle{9\over25}}},46{{\textstyle{25\over36}}}, 68{{\textstyle{32\over49}}},
92{{\textstyle{41\over64}}}\dots
\right\} \hbox{ for } n\in\{1,2,3,\dots\}.
\end{equation}
Consequently, comparing with 
\begin{equation}
\bar p_n \in \left\{2,  2\half, 3\third, 4\quarter, 5{{\textstyle{3\over5}}}, 6{{\textstyle{5\over6}}}, 
8{{\textstyle{2\over7}}}, 9{{\textstyle{5\over8}}}, 11{{\textstyle{1\over9}}}, 
12{{\textstyle{9\over10}}},\dots
\right\} \hbox{ for } n\in\{1,2,3,\dots\},
\end{equation}
we explicitly see that
\begin{equation}
\bar\pi(m^2) \in \{ 0,3,5,7,11, 16,\dots\} \hbox{ for } m\in\{1,2,3,\dots\}.
\end{equation}
That is,  the prime-average version of Legendre actually holds for integer $m\geq 1$. 

Furthermore, the argument above nowhere explicitly uses the fact that $m$ need be integer --- so the prime-average version of Legendre actually holds for real $m\geq \sqrt{5}$.  Finally to check real range $m<\sqrt{5}$ use a suitable finite truncation of the sum
\begin{equation}
\bar\pi(x) = \sum_{i=1}^\infty \Heaviside(x - \bar p_i); \qquad \Heaviside(0)=1.
\end{equation}
Thereby one verifies  that the prime-average version of Legendre certainly holds for real $m\geq 1$. (Actually $m\geq \sqrt{\bar p_2}-1 = \sqrt{5/2}-1= 0.581138830$ is sufficient.)

For another way of extending the argument above, consider this: If we were to define $\bar\pi([m+1]^2)=n_*$, then the average-prime gaps \emph{below} $p_{n_*}$ would be \emph{at most} of size $\ln n_*$ and so
\begin{equation}
\bar\pi([m+1]^2)-\bar\pi(m^2) \geq {2m+1\over \ln n_*} 
= {2m+1\over  \ln \bar \pi([m+1]^2)}  \geq {2m+1\over \ln \{ 2\pi([m+1]^2)\}}.
\end{equation}
Thence for $[m+1]^2 > \exp\left(139/125\right)$, certainly for $m\geq 1$, we have
\begin{equation}
\bar\pi([m+1]^2)-\bar\pi(m^2) 
\geq {2m+1\over \ln \{ 2[m+1]^2/(\ln([m+1]^2)-{139\over125})\}}.
\end{equation}
Expanding the logarithm
\begin{equation}
\bar\pi([m+1]^2)-\bar\pi(m^2) 
\geq {2m+1\over \{ 2\ln[m+1]+\ln 2 -\ln(2\ln[m+1]-{139\over125})\}}.
\end{equation}

Ultimately, for $m\geq 4$ we certainly have
\begin{equation}
\bar\pi([m+1]^2)-\bar\pi(m^2) 
\geq {2m+1\over 2\ln[m+1]}.
\end{equation}
Sacrificing a little precision in the interests of a slightly wider range of validity, for $m\geq 1$ we have
\begin{equation}
\bar\pi([m+1]^2)-\bar\pi(m^2) 
\geq {2m-1\over 2\ln[m+1]}.
\end{equation}

So there will be \emph{many} average-primes $\bar p_n$  between consecutive squares.

\subsection{Prime-average analogue of Oppermann}

The ordinary Oppermann conjecture is the hypothesis that the ordinary primes $p_n$ satisfy $\pi(m[m+1])> \pi(m^2)> \pi(m[m-1])$ for integer $m\geq 2$. This is completely equivalent to demanding that $\pi([m+{1\over2}]^2)> \pi(m^2)> \pi([m-{1\over2}]^2)$.
 
Then the  most compelling prime-average analogue of Oppermann would be that the averaged primes $\bar p_n$ satisfy $\bar\pi([m+{1\over2}]^2)> \bar\pi(m^2) > \bar\pi([m-{1\over2}]^2)$. But (adapting the argument presented above for Legendre) this is easily checked to be true, and in fact much more can be said.

Let $\bar\pi(m^2)=n$, then for integer $n\geq 3$ corresponding to $m\geq \sqrt{5}$ we still have
\begin{equation}
m^2< \bar p_{n+1}   < m^2 + \ln 2 + 2\ln m. 
\end{equation}
Thence
\begin{equation}
m^2< \bar p_{n+1}   <   \left(m+{1\over2}\right)^2 - \left[m-2\ln m+{1\over 4}-\ln 2\right] < \left(m+{1\over2}\right)^2 .
\end{equation}
It is easy to check that for small integers $m$ this continues to hold for $m\geq 1$.
(For real values of $m$ the domain of validity is $m\geq \sqrt{\bar p_1} -{1\over2} = \sqrt{2}-{1\over2} = 0.9142135620$.)

Similarly let $\bar\pi([m-{1\over2}]^2)=n_*$, then through an entirely analogous argument
\begin{equation}
\left(m-{1\over2}\right)^2 < \bar p_{n_*+1}   < \left(m-{1\over2}\right)^2 + \ln 2 + 2\ln\left(m-{1\over2}\right). 
\end{equation}
Thence
\begin{equation}
\left(m-{1\over2}\right)^2 < \bar p_{n_*+1}   <   m^2 - \left[m-2\ln\left(m-{1\over2}\right)+{1\over 4}-\ln 2\right] < m^2 .
\end{equation}
It is easy to check that for small integers $m$ this continues to hold for $m\geq 2$.
(For real values of $m$ the domain of validity is $m\geq \sqrt{\bar p_1}= \sqrt{2} = 1.414213562$.)

So the prime-average analogue of Opperman is unassailably true for real $m\geq 2$.

\subsection{Prime-average analogue of Brocard}

The ordinary Brocard conjecture is the hypothesis that the ordinary primes $p_n$ satisfy $\pi(p_{n+1}^2 )- \pi(p_n^2)\geq  4$ for $n\geq 2$. Note that the ordinary Brocard conjecture is implied by the ordinary Oppermann conjecture.  Note that $p_{n+1}\geq p_n+2$ and that ordinary Oppermann implies the existence of at least one prime in each of the 4 regions 
\begin{equation}
\textstyle{
(p_n^2, [p_n+{1\over2}]^2),  \;\; ([p_n+{1\over2}]^2, [p_n+1]^2),  \;\; 
([p_n+1]^2, [p_n+{3\over2}]^2),  \;\;  ([p_n+{3\over2}]^2, [p_n+2]^2).
}
\end{equation}
Generalizing this, the ordinary Oppermann conjecture implies a generalization of the ordinary Brocard conjecture 
\begin{equation}
\pi(p_{n+1}^2 )- \pi(p_n^2)\geq  2 g_n; \qquad (n \geq 2).
\end{equation}

Then the  most compelling prime-average analogue of Brocard would be that the averaged primes $\bar p_n$ satisfy $\bar\pi(\bar p_{n+1}^2)- \bar\pi(\bar p_n^2) \geq 4$ for suitable values of $n$. But (adapting the argument presented above for Legendre) this is easily checked to be true, and in fact much more can be said.

Write $\bar p_{n+1} = \bar p_n + \bar g_n$ and use prime-average analogue of Oppermann. Then
\begin{equation}
\bar\pi(\bar p_{n+1}^2)- \bar\pi(\bar p_n^2)= \bar\pi( [\bar p_n + \bar g_n]^2)- \bar\pi(\bar p_n^2) > 2\bar g_n;
\qquad 
(n \geq 1). 
\end{equation}
But $\bar g_n > {\ln n\over 2}$ for $n\geq 1$ so 
\begin{equation}
\bar\pi(\bar p_{n+1}^2)- \bar\pi(\bar p_n^2) > \ln n; 
\qquad 
(n \geq 1). 
\end{equation}
So certainly 
\begin{equation}
\bar\pi(\bar p_{n+1}^2)- \bar\pi(\bar p_n^2) > 4; 
\qquad 
(n \geq e^4). 
\end{equation}
Noting that  $e^4 = 54.59815003$, to complete the analysis it suffices to explicitly check the first 55 average primes.
Then
\begin{equation}
\bar\pi(\bar p_{n+1}^2)- \bar\pi(\bar p_n^2) > 4; 
\qquad 
(n \geq 4; \; \bar p_n \geq 4\quarter). 
\end{equation}
So certainly the prime-average analogue of Brocard holds over a suitable range.

\subsection{Prime-average analogue of Firoozbakht}

The ordinary Firoozbakht conjecture is the hypothesis that the quantity $(p_n)^{1/n}$ is monotone decreasing for the ordinary primes $p_n$. The prime-average analogue of  Firoozbakht would be the hypothesis that $(\bar p_n)^{1/n}$ is monotone decreasing for the averaged primes $\bar p_n$. But this is easily checked to be true (see for instance reference~\cite{Sun}) and in fact much more can be said.

\enlargethispage{20pt}
Note
\begin{equation}
(\bar p_{n+1})^{1/(n+1)} < (\bar p_n)^{1/n} \quad \iff \quad  n \ln \bar p_{n+1} < (n+1) \ln \bar p_n
\end{equation}
So let us compute
\begin{equation}
Q:= n \ln \bar p_{n+1}-(n+1)\ln \bar p_n = n \ln (\bar p_{n} +\bar g_n)-(n+1)\ln \bar p_n
= n \ln(1+\bar g_n/\bar p_n) - \ln \bar p_n. 
\end{equation}
This quantity is certainly less than
\begin{equation}
Q < {n \bar g_n\over  \bar p_n} - \ln \bar p_n < {n \ln n \over {1\over2} n\ln n \left[ 1+{\ln\ln n \over\ln n}\right] }- \ln \bar p_n 
< 2 -\ln \bar p_n,
\end{equation}
which clearly becomes negative for $\bar p_n > e^2$, that is, for $n \geq 12$. Checking $n\in[1,11]$ by direct computation shows that $Q<0$ for all $n\geq 1$, so the prime-average analogue of Firoozbakht is unassailably true for all $n\geq 1$. 

\subsection{Prime-average analogue of Fourges}

The ordinary Fourges conjecture is typically presented in terms of first rearranging the ordinary Firoozbakht conjecture into the form
\begin{equation}
\left( \ln p_{n+1} \over \ln p_n \right)^n < \left(1+{1\over n}\right)^n;\qquad (n \geq 1),
\end{equation}
and then \emph{weakening} it by making the \emph{less restrictive} demand that
\begin{equation}
\left( \ln p_{n+1} \over \ln p_n \right)^n <  e;\qquad (n \geq 1),
\end{equation}

To obtain a prime-average analogue of the Fourges conjecture we first rearrange the prime-average analogue Firoozbakht conjecture into the form
\begin{equation}
\left( \ln \bar p_{n+1} \over \ln \bar p_n \right)^n < \left(1+{1\over n}\right)^n;\qquad (n \geq 1),
\end{equation}
and then \emph{weaken} it by making the \emph{less restrictive} demand that
\begin{equation}
\left( \ln \bar p_{n+1} \over \ln \bar p_n \right)^n <  e;\qquad (n \geq 1).
\end{equation}
Since this is a weakening, the prime-average analogue of Fourges is unassailably true for all $n\geq 1$. 

\subsection{Prime-average analogue of Nicholson}

The ordinary Nicholson conjecture is typically presented in terms of first rearranging the ordinary Firoozbakht conjecture into the form 
\begin{equation}
\left( p_{n+1} \over p_n \right)^n <  p_n;\qquad (n \geq 1),
\end{equation}
and then (slightly) \emph{strengthening} it by making the \emph{more restrictive} demand that
\begin{equation}
\left( p_{n+1} \over p_n \right)^n <  n \ln n;\qquad (n \geq 1).
\end{equation}

To obtain a prime-average analogue of the Nicholson conjecture we first rearrange the prime-average analogue Firoozbakht conjecture into the form
\begin{equation}
\left( \bar p_{n+1} \over \bar p_n \right)^n <  \bar p_n;\qquad (n \geq 1),
\end{equation}
and then, noting $\bar p_n > p_{[n/2]} > (n/2)\ln(n/2)$, with the final inequality holding for $n\geq 20$,  (slightly) \emph{strengthen} it by making the \emph{more  restrictive} demand that
\begin{equation}
\left( \bar p_{n+1} \over \bar p_n \right)^n <  (n/2)\ln (n/2) ;\qquad (n \geq n_0),
\end{equation}
for some $n_0\leq 20$ yet to be determined.
That is, our proposed form of the prime-average analogue of Nicholson is equivalent to 
\begin{equation}
n \ln \left( \bar p_{n+1} \over \bar p_n \right) <  \ln\left( (n/2) \ln (n/2)\right) ;\qquad (n \geq n_0).
\end{equation}

Since
\begin{equation}
n \ln \left( \bar p_{n+1} \over \bar p_n \right) = n \ln \left(1+ {\bar g_{n} \over \bar p_n} \right)
<  n \ln \left(1+{2\over n}\right) < 2,
\end{equation}
with the inequalities holding for $n\geq 4$, we see that our proposed version of the prime-average analogue of the Nicholson conjecture will plausibly hold whenever
\begin{equation}
\ln\left( (n/2) \ln (n/2)\right) > 2,
 \end{equation}
 that is, whenever $n \geq 5$, and so will certainly hold for all $n\geq 20$. Finally, explicitly checking all integers below 20 we have an explicit  prime-average analogue of Nicholson:
 \begin{equation}
\left( \bar p_{n+1} \over \bar p_n \right)^n <  (n/2)\ln (n/2) ;\qquad (n \geq 6).
\end{equation}

\subsection{Prime-average analogue of Farhadian}

\enlargethispage{20pt}
The ordinary Farhadian conjecture is typically presented in terms of  (slightly) \emph{strengthening} the ordinary Nicholson conjecture by making the even \emph{more restrictive} demand that
\begin{equation}
\left( p_{n+1} \over p_n \right)^n <  {p_n \ln n\over\ln p_n};\qquad (n \geq 1).
\end{equation}

To obtain a prime-average analogue of the Farhadian conjecture we again  (slightly) \emph{strengthen} the prime-average analogue of Nicholson conjecture by making  \emph{some} (slightly) more restrictive demand.  There is a potential infinity of stronger demands that one might make, but we shall try to keep close to the spirit of the original Farhadian conjecture by demanding
\begin{equation}
\left( \bar p_{n+1} \over \bar p_n \right)^n <  {p_n\over 2 \ln p_n} \: \ln (n/2) ;\qquad (n \geq n_0).
\end{equation}
Since $p_n/\ln p_n <n$ for $n\geq 4$ this certainly strengthens the prime-average analogue of Nicholson. 

This would be equivalent to
\begin{equation}
n \ln\left( \bar p_{n+1} \over \bar p_n \right) <  \ln\left({p_n\over 2 \ln p_n} \: \ln (n/2) \right);\qquad (n \geq n_0).
\end{equation}

As previously noted 
\begin{equation}
n \ln \left( \bar p_{n+1} \over \bar p_n \right) = n \ln \left(1+ {\bar g_{n} \over \bar p_n} \right)
<  n \ln \left(1+{2\over n}\right) < 2,
\end{equation}
with these inequalities holding for $n\geq 4$. 

\enlargethispage{40pt}
This now suggests we  consider the inequality
\begin{equation}
\ln\left({p_n\over 2 \ln p_n} \: \ln (n/2) \right) > 2
\end{equation}
which holds for $n\geq 5$, (though the original weakening to get to get to prime-average analogue Nicholson holds only for $n\geq 20$). Explicitly checking small integers
\begin{equation}
n \ln\left( \bar p_{n+1} \over \bar p_n \right) <  \ln\left({p_n\over 2 \ln p_n} \: \ln (n/2) \right);\qquad (n \geq 7).
\end{equation}
Equivalently 
\begin{equation}
\left( \bar p_{n+1} \over \bar p_n \right)^n <  {p_n\over 2 \ln p_n} \: \ln (n/2) ;\qquad (n \geq 7).
\end{equation}
As expected, because the prime-average analogue of Farhadian is slightly stronger than the prime-average analogue of Nicholson, it is valid only on a slightly smaller domain.

\section{Discussion}\label{S:discussion}

We have explicitly seen above that the averaged primes $\bar p_n ={1\over n}  \sum_{i=1}^n p_i$ satisfy \break suitably defined averaged-prime analogues of 
the Cramer, Andrica, Legendre, \break Oppermann, Fourges, Firoozbakht, Nicholson, and Farhadian hypotheses. The key step in all cases was the extremely tight bound in the averaged-prime gaps: $\bar g_n = \bar p_{n+1} - \bar p_n < \ln n$ for $n\geq 3$. 
The only (relatively minor) potential difficulty lies in the design and formulation of truly compelling prime-average analogues of these hypotheses.

On the other hand, what does this tell us about the ordinary primes $p_n$?\\
Not as much as one might hope. We certainly have the identity
\begin{equation}
g_n  = [(n+1)\bar p_{n+1}- n\bar p_n] - [n \bar p_{n}- (n-1)\bar p_{n-1}] 
= [\bar g_n +\bar g_{n-1}]+  n [\bar g_n -\bar g_{n-1}].
\end{equation}
Unfortunately a naive application of this equality merely leads to $g_n = \O(n\ln n)$ \break --- bounds too weak to be useful. So even though the numerical evidence is quite compelling~\cite{Verifying1,Verifying2}, we must conclude that significant progress on any of the usual Cramer, Andrica, Legendre, Oppermann, Fourges, Firoozbakht, Nicholson, and \break Farhadian conjectures for ordinary primes would require considerably more subtle arguments. 

\bigskip

\hrule\hrule\hrule

\clearpage

\end{document}